\documentclass[11pt] {article} 

  \usepackage{amsmath}
    \usepackage{amssymb}
    
    \usepackage{pdfsync}

\newtheorem{example}{Example}[section]
\newtheorem{theorem}{Theorem}[section]
\newtheorem{lemma}{Lemma}[section]
\newtheorem{corollary}{Corollary}[section]

\newtheorem{remark}{Remark}[section]
\newcommand{\eqnsection}{
   \renewcommand{\theequation}{\thesection.\arabic{equation}}
   \makeatletter
   \csname @addtoreset\endcsname{equation}{section}
   \makeatother}
 
\def \ov{\overline}

\def \be{\begin{equation}}
\def \ee{\end{equation}}
\def \bt{\begin{theorem}} 
\def \et{\end{theorem}}
\def \bl{\begin{lemma}} 
\def \el{\end{lemma}}
\def \bea{\begin{eqnarray}}
\def \eea{\end{eqnarray}}
\def \bas{\begin{eqnarray*}}
\def \eas{\end{eqnarray*}}

\def \al{\alpha}
\def \bb{\beta}
\def \ga{\gamma}
\def \Ga{\Gamma}

\def \ff{\infty}
\def \wh{\widehat}
\def \wt{\widetilde}

\def \cd{\,\cdot\,}

\def \BB{{\mathcal B}}

\def\stl{\stackrel{law}{=}}
\def \({\left(}
\def \){\right)}

\def \nn{\nonumber}
\def \Proof{\noindent{\bf Proof $\,$ }}

\def \bc{\begin{center} }
\def \ec{\end{center} }
\def \bs{\begin{slide} }
\def \es{\end{slide} }

\def\square{{\vcenter{\vbox{\hrule height.3pt
        \hbox{\vrule width.3pt height5pt \kern5pt
           \vrule width.3pt}
        \hrule height.3pt}}}}
\def\qed{{\hfill $\square$ \bigskip}}

\eqnsection
\begin{document}

\title{ Permanental random variables, ${M}$-matrices and  $\al$-permanents}

  \author{  Michael B. Marcus\,\, \,\, Jay Rosen \thanks{Research of     Jay Rosen was partially supported by  grants from the National Science Foundation.   }}
\maketitle

\begin{abstract}   We explore some properties of a recent representation  of permanental vectors which expresses them as sums of independent  vectors with components that are independent gamma random variables.\end{abstract}

 \maketitle

\section{Introduction}\label{sec-1} 

An $R^{n}$ valued   $\al$-permanental random variable $X=(X_{1},\ldots, X_{n})$ is  a random variable with Laplace transform 
\begin{equation}
   E\(e^{-\sum_{i=1}^{n}s_{i}X_{i}}\) 
 = \frac{1}{ |I+RS|^{ \al}}   \label{int.1} 
 \end{equation}
for some $n\times n$ matrix $R$   and diagonal matrix $S$ with entries $s_{i}$, $1\le i\le n$, and $\al>0$. Permanental random variables were introduced by Vere-Jones, \cite{VJ}, who called them multivariate gamma distributions. (Actually he considered the moment generating function.)

   In \cite[Lemma 2.1]{MRcond} we obtain a representation for permanental random variables with the property that $A=R^{-1}$ is an $M$-matrix.
A matrix
$A=\{ a_{ i,j}\}_{ 1\leq i,j\leq n}$  is said to be a nonsingular
$M$-matrix   if
\begin{enumerate}
\item[(1)] $a_{ i,j}\leq 0$ for all $i\neq j$.
\item[(2)] $A$ is nonsingular and $A^{ -1}\geq 0$.
\end{enumerate}

The representation depends on the $\al$-perminant of the off diagonal elements of $A$ which we now define.

  The $\al$-perminant of $n\times n$ matrix $M$ is
\begin{equation}
|M|_{\al}= \begin{vmatrix}
M_{1,1} & \cdots& M_{1,n}\\ 
\cdots& &\cdots \\
M_{n,1} &\cdots&M_{n,n}
\end{vmatrix}_{\al}=\sum_{\pi}\al^{c(\pi)}M_{1,\pi(1)}M_{2,\pi(1)}\cdots M_{n,\pi(n)}.\label{vj.71}
\end{equation}
Here the sum runs over all permutations $\pi$ on $[1,n]$ and $c(\pi)$ is the number of cycles in $\pi$.

 We use boldface, such as  ${\bf x }$, to denote vectors. Let   ${\bold k}=(k_{1},\ldots, k_{n})\in \mathbb{N}^{n}$ and $ |\bold k|=\sum_{l=1}^{n}k_{l}$. For  $1\leq m\leq  |\bold k|$,   set    $ i_{m}=j$, where 
 \begin{equation}
\sum_{l=1}^{j-1}k_{l}< m\leq \sum_{l=1}^{j}k_{l}.\label{vj.70v}
 \end{equation}
For any $n\times n$ matrix  $C=\{ c_{ i,j}\}_{ 1\leq i,j\leq n}$  we define  
  \begin{equation}
C(\bold k)= \begin{bmatrix}
c_{ i_{1}, i_{1}} & c_{ i_{1}, i_{2}}&\cdots &c_{i_{1}, i_{|\bold k|}}\\ 
c_{ i_{2}, i_{1}} & c_{ i_{2}, i_{2}}&\cdots &c_{i_{2}, i_{|\bold k|}}\\ 
\cdots& &\cdots \\
c_{ i_{|\bold k|}, i_{1}} & c_{ i_{|\bold k|}, i_{2}}&\cdots &c_{i_{|\bold k|}, i_{|\bold k|}}
\end{bmatrix},\label{vj.71va}
\end{equation}
 and $C(\bold 0)=1$.	For example,  if ${\bold k}=(0,2, 3)$ then $|\bold k|=5$ and $i_{1}=i_{2}=2$ and $i_{3}=i_{4}=i_{5}=3$, 
   \begin{equation}
C(0,2, 3)= \begin{bmatrix}
c_{ 2, 2} & c_{2, 2}&c_{2, 3}&c_{2, 3}&c_{2, 3}\\ 
c_{ 2, 2} & c_{2, 2}&c_{2, 3}&c_{2, 3}&c_{2, 3}\\  
c_{ 3, 2} & c_{3, 2}&c_{3, 3}&c_{3, 3}&c_{3, 3}\\ 
c_{ 3, 2} & c_{3, 2}&c_{3, 3}&c_{3, 3}&c_{3, 3}\\ 
c_{ 3, 2} & c_{3, 2}&c_{3, 3}&c_{3, 3}&c_{3, 3}
\end{bmatrix}.\label{vj.71vs}
\end{equation}

 Here is an alternate description of $C(\bold k)$.  For any $n\times n$ matrix  $C=\{ c_{ i,j}\}_{ 1\leq i,j\leq n}$ the matrix $C(\bold k)$ is an $|\bold k|\times |\bold k|$ matrix with its first $k_{1}$ diagonal elements equal to $c_{1,1}$, its next $k_{2}$ diagonal elements equal to $ c_{2,2}$, and so on. The general element   $C(\bold k)_{  p, q }= c_{ \bar p, \bar q }$, where $\bar   p$ is equal to either  index of    diagonal element in row $p$, (the diagonal element has  two  indices but they are the same), and $\bar   q$ equal to either index of the diagonal element in column $q$. Thus in the above example we see that   $C(0,2,3)_{ 4, 1 }= c_{3, 2}$.

\medskip	 Suppose that $A$ is an $n\times n$ $M$-matrix. Set  $a_{i}=a_{i,i}$ and  write
    \begin{equation}
A=D_{A}-B,\label{vj.73}
\end{equation}
 where $D$ is a diagonal matrix with entries  $a_{1},\ldots, a_{n} $   and  all the elements of $B$ are non-negative. (Note that   all the diagonal elements of $B$ are equal to zero). In addition set
 \begin{equation}
\ov A=D_{A}^{-1}A =I-D_{A}^{-1}B:=I-\ov B.\label{vj.73q}
\end{equation}
 	
The next lemma is   \cite[Lemma 2.1]{MRcond}. 
\begin{lemma} \label{lem-1.1}
Let $A=R^{-1}$ be an $n\times n$ nonsingular  M-matrix with diagonal entries $a_{1},\ldots, a_{n}$ and   $S$ be an $n\times n$ diagonal matrix with entries $(s_{1},\ldots, s_{n})$. Then (\ref{int.1}) is equal to
\bea
\lefteqn{
\frac{|A|^{\al}}{\prod_{i=1}^{n}a^{\al}_{i}}\sum_{{\bold k}=(k_{1},\ldots, k_{n}) }  {|B(\bold k)|_{\al}\over  \prod_{i=1}^{n}a_{i}^{k_{i}} k_{i}! }  \, \frac{1}{ (1+(s_{1}/a_{1}))^{\al+k_{1}}\cdots  (1+(s_{n}/a_{n}))^{\al+k_{n}}} \nn}\\
 &&=  
 |\ov A|^{\al}\sum_{{\bold k}=(k_{1},\ldots, k_{n}) } \frac{   |\ov B(\bold k)|_{\al}}{\prod_{i=1}^{n}k_{i}!} \, \frac{1}{ (1+(s_{1}/a_{1}))^{\al+k_{1}}\cdots  (1+(s_{n}/a_{n}))^{\al+k_{n}}}.\label{1.11}
 \eea 
where   the sum is over all ${\bold k}=(k_{1},\ldots, k_{n})\in \mathbb{N}^{n}$.   (The series converges for all $s_{1},\ldots,s_{n}\in R_{+}^{n}$ for all $\al>0$.)
  \end{lemma}
  
  Setting $S=0$ we see that 
  \begin{equation}
  \frac{ |  A|^{\al}}{\prod_{i=1}^{n}a_{i}^{\al}}\sum_{{\bold k}=(k_{1},\ldots, k_{n}) } \frac{   |  B(\bold k)|_{\al}}{\prod_{i=1}^{n}a_{i}^{\al}k_{i}!}=   |\ov A|^{\al}\sum_{{\bold k}=(k_{1},\ldots, k_{n}) } \frac{   |\ov B(\bold k)|_{\al}}{\prod_{i=1}^{n}k_{i}!}=1.\label{1.9}
   \end{equation}
Let   $Z_{\al,\ov B} $ be an  $n$-dimensional integer valued random variable with
  \begin{equation}
  P\(Z_{\al,\ov B}=(k_{1},\ldots, k_{n})\)=    |\ov A|^{\al}  \frac{   |\ov B(\bold k)|_{\al}}{\prod_{i=1}^{n}k_{i}!}.\label{9.11w}
  \end{equation} 
  (We  omit writing the subscript $ {\al,\ov B}$ when they are fixed from term to term.)

  The sum in (\ref{1.11}) is the Laplace transform of the the $\al$-permanental random variable $X$. Therefore, 
     \bea  X&\stl&\sum_{{\bold k}=(k_{1},\ldots, k_{n}) } I_{   k_{1},\ldots, k_{n}  } (Z)\(\xi_{ \al+k_{1},a_{1}},\ldots,\xi_{\al+k_{n},a_{n}}\) \label{9.10}\\
&\stl&  \(\xi_{ \al+Z_{1},a_{1}},\ldots,\xi_{\al+Z_{n},a_{n}}\),\nn
  \eea  
  where $Z$ and all the gamma distributed random variables, $\xi_{u,v}$  are independent.   Recall that the probability density function of $\xi_{u,v}$ is
	 \begin{equation}
 f(u,v;x)  =\frac{v^{u}x^{u-1}e^{-vx}}{\Ga(u)}\quad\mbox{for $x> 0$ and $u,v>0$},
   \end{equation}
 and equal to 0 for   $x\le 0$. 
 We see that when $X$ has probability density $\xi_{u,v}$,  $vX$ has probability density $\xi_{u,1}$.
It is easy to see that
 \begin{equation}
   E(\xi_{u,1}^{p})=\frac{\Ga(p+u)}{\Ga( u)}.\label{1.13}
   \end{equation}

 It follows from (\ref{9.10}) that
  for measurable functions  $f$ on  $R^{n}_{+}$,
   \bea  E (f(X))\nn 
 &=  &  \sum_{{\bold k}=(k_{1},\ldots, k_{n}) }  P\(Z=(k_{1},\ldots, k_{n})\)  E \( f\(\xi_{\al+k_{1},a_{1}} ,\ldots,\xi_{\al+k_{n},a_{n}} \)\)\\
 &=&E \( f\(\xi_{\al+Z_{1},a_{1}} ,\ldots,\xi_{\al+Z_{n},a_{n}} \)\).\label{9.12} 
  \eea 
 
 Since
 \begin{equation}
   \xi_{\al+\bb,a}\stl    \xi_{\al,a}  + \xi_{ \bb,a},
   \ee
it  follows from (\ref{9.12}) that for all increasing functions $f$
 \begin{equation}
   E (f(X))\ge   E \( f\(\xi_{\al,a_{1}} ,\ldots,\xi_{\al ,a_{n}} \)\) .\label{1.14}
   \end{equation}
  We explain in  \cite{MRcond} that in some respects  (\ref{1.14}) is a generalization of the Sudakov Inequality for Gaussian processes and use it to obtain sufficient conditions for permanental processes to be unbounded. 
  
 A permanental process is a process with finite joint distributions that are permanental random variables. For example, let $G=\{G(t),t\in R\}$ be a Gaussian process with covariance $\wt R(s,t)$. Then for all $n$ and all $t_{1}, \ldots,t_{n}$ in $R^{n}$,  $(G^{2}(t_{1})/2, \ldots, G^{2}(t_{n})/2)$ is an $n$-dimensional 1/2-permanental random variable, with $R$  in (\ref{int.1}) equal to    the kernel	$\{\wt R(t_{i},t_{j})\}_{i,j=1}^{n}$.  The stochastic process $G^{2}/2=\{G^{2}(t)/2,t\in R\}$ is a 1/2-permanental process. In \cite{MRcond} we consider permanental processes defined for all $\al>0$ and for kernels $R(s,t)$ that need not be symmetric. 
  
  \medskip	In the first part of this paper we give some properties of   the random variable $Z$. It turns out that it is easy to obtain the Laplace transform of $Z$.

\begin{lemma} 
\begin{equation}
   E\(e^{-\sum_{i=1}^{n}s_{i}Z_{i}}\) 
 = \frac{|\ov A|^{\al}}{ |I-(\ov B E({\bf s}))|^{ \al}}   \label{int.11} 
 \end{equation}
where $E({\bf s})$ is an $n\times n$ diagonal matrix with entries $e^{-s_{i}}$, $i=1,\ldots,n$.
 \end{lemma}
  
  \Proof 
  \bea
    E\(e^{-\sum_{i=1}^{n}s_{i}Z_{i}}\)& =& \sum_{{\bold k}=(k_{1},\ldots, k_{n}) }  e^{-\sum_{i=1}^{n}s_{i}k_{i}}  P\(Z=(k_{1},\ldots, k_{n})\)\nn\\
    & =& |\ov A|^{\al}  \sum_{{\bold k}=(k_{1},\ldots, k_{n}) } \prod_{i=1}^{n} e^{-s_{i}k_{i}}    \frac{   |\ov B(\bold k)|_{\al}}{\prod_{i=1}^{n}k_{i}!}.\label{1.5}
   \eea
  Note that 
  \begin{equation}
    \prod_{i=1}^{n} e^{-s_{i}k_{i}}    \frac{   |\ov B(\bold k)|_{\al}}{\prod_{i=1}^{n}k_{i}!}=    \frac{   |(\ov {B} E ({\bf s}))(\bold k)|_{\al}}{\prod_{i=1}^{n}k_{i}!}.
   \end{equation}
  By (\ref{1.9})   with $\ov B(0)$ replaced by $\ov BE(\bf s)$ for each fixed ${\bf s}$
  \begin{equation}
|I-\ov {B} E({\bold s})|^{\al}    \sum_{{\bold k}=(k_{1},\ldots, k_{n}) }       \frac{   | (\ov {B} E ({\bf s}))(\bold k)|_{\al}}{\prod_{i=1}^{n}k_{i}!}=1.
   \end{equation}
We get  (\ref{int.11})   from this and   (\ref{1.5}).\qed

\medskip	A significant property of permanental random variables  that are determined by kernels that are inverse $M$-matrices is that they are infinitely divisible. Similarly it follows from (\ref{int.11}) that for all $\al,\bb>0$
\begin{equation}
   Z_{\al+\bb, \ov B}\stl  Z_{\al , \ov B} +Z_{ \bb, \ov B}.
   \end{equation}

\medskip	We can differentiate (\ref{int.11}) to give a simple formula for the moments of the components of $Z_{\al,\ov B}$, which we simply denote by $Z$ in the following lemma.

  \begin{lemma}\label{lem-1.3} For any   integer $m\ge 1$ and $1\le p\le n$,  
\be E(Z_{p}^{m})= \label{1.17}  \sum_{\stackrel{ j_{0}+\cdots j_{l}=m, \,j_{i}\ge 1}  {l=0,1,\ldots, m-1} }(-1)^{l+m+1}\al^{j_{0}}(\al+1) ^{j_{1}}\cdots (\al+l) ^{j_{l }}(R_{p,p}A_{ p,p})^{\,l}(  R_{p,p}A_{p,p}-1).
 \ee 
   ($A_{p,p}$ is also referred to as $a_{p}$ elsewhere in this paper.)
 \end{lemma}

 	  \Proof  To simplify the notation we take $p=1$. Note that by (\ref{int.11}) and the fact that $\ov A=I-(\ov B E({\bf 0}))$
	  \bea
  E(Z_{1}^{m})
& =& (-1)^{m}     \frac{\partial^{m}}{\partial s^{m}_{1}}  \(\frac{|\ov A|^{\al}}{ |I-(\ov B E({\bf s}))|^{ \al}} \) \Bigg|_{{\bf s}=0} \label{int.11j} \\
  & =& (-1)^{m}  |I-(\ov B E({\bf 0}))|^{ \al}   \frac{\partial^{m}}{\partial s^{m}_{1}}  \(\frac{1}{ |I-(\ov B E({\bf s}))|^{ \al}} \)\Bigg|_{{\bf s} =0}.\nn
 \eea
Hence to prove  (\ref{1.17}) it suffices to show that   for any $m$
   \begin{eqnarray}
 \lefteqn{ \frac{\partial^{m}}{\partial s^{m}_{1}}\( \frac{1}{ |I-(\ov B E({\bf s}))|^{ \al}}\)\Bigg|_{{\bf s} =0}=\sum_{\stackrel{ j_{0}+\cdots j_{l}=m, \,j_{i}\ge 1}  {l=0,1,\ldots,m-1} }(-1)^{l+1}
   \label{ind3x}}\\
   &&\hspace{.5in}   \al^{j_{0}}(\al+1) ^{j_{1}}\cdots (\al+l) ^{j_{l }}\frac{(R_{1,1}A_{1,1})^{\,l}}{ |I-(\ov B E({\bf 0}))|^{\al}}(  R_{1,1}A_{1,1}-1).\nonumber
   \end{eqnarray}

Note that for any  $\ga>0$ we have
  \begin{equation}
  \frac{\partial}{\partial s_{1}}\( \frac{1}{ |I-(\ov B E({\bf s}))|^{ \ga}}\)= - \frac{\ga}{ |I-(\ov B E({\bf s}))|^{ \ga+1}} \frac{\partial}{\partial s_{1}}  |I-(\ov B E({\bf s}))| .\label{1.18}
   \end{equation}
  We expand the determinant by the first column. Since   $ \ov b_{1,1}=0$ we have  
   \begin{equation}
    |I-(\ov B E({\bf s}))| = M_{1,1}- \ov b_{2,1}e^{-s_{1}}M_{2,1}+\ov b_{3,1}e^{-s_{1}}M_{3,1} \cdots \pm \ov b_{n,1}e^{-s_{1}}M_{n,1}\label{1.19}
   \end{equation}
  where $M_{i,1}$ are minors of $I-(\ov B E({\bf s}))$ and the last sign is plus or minus according to whether $n$ is odd or even. Note that the terms $M_{i,1}$ are not functions of $s_{1}$. 
 Using  (\ref{1.19}) we see that
  \bea
   \frac{\partial}{\partial s_{1}}  |I-(\ov B E({\bf s}))|&= &   \ov b_{2,1}e^{-s_{1}}M_{2,1}-\ov b_{3,1}e^{-s_{1}}M_{3,1} \cdots \mp \ov b_{n,1}e^{-s_{1}}M_{n,1}\nn\\
   &=&-   |I-(\ov B E({\bf s}))|+M_{1,1}. \label{jmo}
   \eea
 Using this    we get 
   \bea
  \frac{\partial}{\partial s_{1}}\( \frac{1}{ |I-(\ov B E({\bf s}))|^{ \ga}}\)&= & \frac{\ga}{ |I-(\ov B E({\bf s}))|^{ \ga }}- \frac{\ga M_{1,1} }{ |I-(\ov B E({\bf s}))|^{ \ga+1 }}  \label{1.21}\\
  &= & -\frac{\ga}{ |I-(\ov B E({\bf s}))|^{ \ga }}\( \frac{ M_{1,1} }{ |I-(\ov B E({\bf s}))| }-1\).\nn
   \eea
  
  We now show by induction on $m$ that
    \begin{eqnarray}
\lefteqn{ \frac{\partial^{m}}{\partial s^{m}_{1}}\( \frac{1}{ |I-(\ov B E({\bf s}))|^{ \al}}\)=\sum_{\stackrel{ j_{0}+\cdots j_{l}=m, \,j_{i}\ge 1}  {l=0,1,\ldots,m-1} }(-1)^{l+1}
   \label{ind3}}\\
   &&   \al^{j_{0}}(\al+1) ^{j_{1}}\cdots (\al+l) ^{j_{l }}\frac{ M_{1,1}^{\,l}}{ |I-(\ov B E({\bf s}))|^{\al+ l}}\( \frac{ M_{1,1}}{ |I-(\ov B E({\bf s}))| } -1\).\nonumber
   \end{eqnarray}  
     It is easy to see that  for $m=1$  this agrees with (\ref{1.21}) for $\ga=\al$. Assume that (\ref{ind3}) holds for $m$. We show it holds for $m+1$. We take another derivative with respect to $s_{1}$. It follows from  (\ref{1.21}) that
     \begin{eqnarray}
     &&\frac{\partial}{\partial s_{1}}\(  \frac{ M_{1,1}^{\,l}}{ |I-(\ov B E({\bf s}))|^{\al+ l}}\( \frac{ M_{1,1}}{ |I-(\ov B E({\bf s}))| } -1\)  \)
     \label{ind4}\\
     &&\qquad =\frac{\partial}{\partial s_{1}}\(  \frac{ M_{1,1}^{\,l+1}}{ |I-(\ov B E({\bf s}))|^{\al+ l+1}}- \frac{ M_{1,1}^{\,l}}{ |I-(\ov B E({\bf s}))|^{\al+ l}} \)  \nonumber\\
     &&\qquad =  -\frac{(\al+ l+1) M_{1,1}^{\,l+1}}{ |I-(\ov B E({\bf s}))|^{\al+ l+1}} \(\frac{ M_{1,1} }{ |I-(\ov B E({\bf s}))| }-1\)  \nonumber\\
     &&\qquad\qquad \qquad+\frac{(\al+ l) M_{1,1}^{\,l}}{ |I-(\ov B E({\bf s}))|^{\al+ l}} \(\frac{ M_{1,1} }{ |I-(\ov B E({\bf s}))| }-1\).  \nonumber
     \end{eqnarray}
     
     Let us consider the term corresponding $l=k$   when we take another derivative with respect to $s_{1}$.   Two sets of terms in   (\ref{ind3}) contribute to this. One set are the terms in which $j_{0}+\cdots+j_{k}=m-1$, $j_{k}\ge 1$ which become terms in which $j_{0}+\cdots+(j_{k}+1)=m$, $j_{k}\ge 1$, when 
  $\frac{ M_{1,1}^{\,l}}{ |I-(\ov B E({\bf s}))|^{\al+ l}}\( \frac{ M_{1,1}}{ |I-(\ov B E({\bf s}))| } -1\)$   is replaced by the last line of (\ref{ind4}). This  almost gives us all we need.  We are only lacking $j_{0}+\cdots+j_{k}=m$, $j_{k}= 1$. This comes from the next to last line of (\ref{ind4}) multiplying the terms in (\ref{ind3}) in which  $l=k-1$. One can check that the sign of the terms for $l=k$ for $m+1$ is different from the sign of the  terms for $l=k$ for $m$ which is what we need. This completes the proof by induction.
 
Recall that $\ov A=I-(\ov B E({\bf 0}))$ and that $M_{1,1}$ is actually a function of $\bf s$.
Therefore,
  \begin{equation}
   \frac{ g(\bf 0)} { |I-(\ov B E({\bf 0}))| }=\frac{M_{1,1}(\bf 0)}{|\ov A|}={((\ov A)^{-1})_{1,1}}=R_{1,1}A_{1,1},\label{1.29}
   \end{equation}
  by (\ref{vj.73q}). Combining this with (\ref{ind3})    we get (\ref{ind3x}).
\qed

   Recall that  
  \begin{equation}
   \ov R_{p,p}= R_{p,p}A_{p,p}.
   \end{equation}

  The next lemma gives  relationship between the moments of the components of  a permanental random variables $X$ and the components of the corresponding random variables  $Z$.

\begin{lemma} For $m_{j}\ge 1$,
	We have \bea
 E\(\prod_{j=1}^{n}(a_{j}X_{j})^{m_{j}}\)\label{mom.2} 
=E\(\prod_{j=1}^{n}\prod_{ l=0}^{m_{j}-1}\( {\al+Z_{j}+l} \)\).   \label{1.24}
 \eea
 or, equivalently
 \begin{equation}
    |\ov R(\bold m)|_{\al}=E\(\prod_{j=1}^{n}  \prod_{ l=0}^{m_{j}-1}\( {\al+Z_{j}+l} \) \).\label{2.10}
   \end{equation}
 
\end{lemma}

 \Proof     
  Let $a_{1},\ldots,a_{n}$ denote the diagonal elements of $A$ and set $Y=(a_{1}X_{1},\ldots,\newline a_{n}X_{n})$. Then 
     \be 
  Y\stl \(\xi_{ \al+Z_{1},1},\ldots,\xi_{\al+Z_{n},1}\)\label{9.10w} 
    \ee 
 The left-hand side of (\ref{1.24}) is $ E\(\prod_{j=1}^{n}( Y_{j})^{m_{j}}\)$. Therefore, by (\ref{9.12}) it is equal to
 \be     E\(\prod_{j=1}^{n}\xi_{ \al+Z_{j},1}^{m_{j}} \)  =E\(\prod_{j=1}^{n} {\Ga( \al+Z_{j}+m_{j}) \over \Ga( \al+Z_{j})}\), \ee 
from which we get (\ref{1.24}). 

 		  It follows from   \cite[Prop. 4.2]{VJ}   that for any ${\bf m}=(m_{1},\ldots, m_{n})$ \begin{equation}
E\(\prod_{j=1}^{n}X_{j}^{m_{j}}\)=|R(\bold m)|_{\al}.\label{mom.1}
\end{equation}
Since $|R(\bold m)|_{\al}   \prod_{j=1}^{n}a_{j}^{m_{j}}=    |\ov R(\bold m)|_{\al}$ we get (\ref{2.10}).\qed

One can use the approach of Lemma \ref{lem-1.3} or try to invert (\ref{1.24}) to find mixed moments of $Z_{i}$. Either approach seems difficult. However, it is easy to make a little progress in this direction.

 \begin{lemma} For all $i$ and $j$, including $i=j$,

 \be
 \mbox{Cov } Z_{i} Z_{j}= \mbox{Cov } a_{i }X_{i}a_{j} X_{j}=\al \,a_{1}a_{2}R_{i,j}R_{j,i}\label{1.36}.
 \ee
 \end{lemma}
\Proof By Lemma 1.4 
\begin{equation}
    E(Z_{i})+\al  =\al a_{i}R_{i,i}\label{1.36a}
   \end{equation}and \be
E(a_{i}a_{j}X_{i }X_{j})=E((\al + Z_{i})(\al+ Z_{j}))\label{1.37}.
\ee
We write 
\bea
   E((\al + Z_{i})(\al+ Z_{j}))&=&\al^{2}+\al E( Z_{i})+\al E( Z_{j})+E( Z_{i})E( Z_{j})\label{1.39}\\
   &=&\al^{2}+\al E( Z_{i})+\al E( Z_{j})+Cov( Z_{i}  Z_{j})+E( Z_{i})E(  Z_{j})\nn\\
      &=& (E( Z_{i})+\al) (E( Z_{j})+\al)+Cov( Z_{i}  Z_{j}) \nn\\
      &=& (\al a_{i}R_{i,i}) (\al a_{j}R_{j,j})+Cov( Z_{i}  Z_{j} ) , \nn
   \eea
where we use (\ref{1.36a}) for the last line. Using (\ref{1.39}) and calculating the left-hand side of (\ref{1.37}) we get the equality of the first and third terms in (\ref{1.36}). To find the equality of the second and third terms in (\ref{1.36}) we differentiate the Laplace transform of $(X_{1},X_{2})$ in (\ref{int.1}).\qed

   If $\bold m=(0,\ldots,m_{j}, 0,\ldots 0):=\wt {\bold m}$ it follows from (\ref{2.10}) that
\begin{equation}
   |\ov R(\wt{\bold m})|_{\al}=E\(   \prod_{ l=0}^{m_{j}-1}\( {\al+Z_{j}+l} \) \).\label{2.13}
   \end{equation}
   Note that 
   \begin{equation}
      |\ov R(\wt{\bold m})|_{\al}=\ov R_{j,j} ^{m_{j}}|E_{m_{j}}|_{\al}
   \end{equation}
where $E_{m_{j}}$ is an $m_{j}\times m_{j}$ matrix with all entries equal to 1.  Therefore, by   \cite[Proposition 3.6]{VJ}
\be  
    |\ov R(\wt{\bold m})|_{\al}= \ov R_{j,j} ^{m_{j}}\prod_{l=0}^{m_{j}-1}(\al+l ).\label{2.14}
   \ee

Combining (\ref{2.13})   and (\ref{2.14}) we get the following inversion of (\ref{1.17}):
\begin{lemma} \label{lem-1.6}
\bea
    \ov R_{j,j} ^{m_{j}}  \label{mom.22} 
=E\( \prod_{ l=0}^{m_{j}-1}\frac{\al+Z_{j}+l}{\al+l} \) .   
 \eea 
  \end{lemma}

As a simple example of  (\ref{mom.22}) or (\ref{1.24}) we have 
 \begin{equation}
 \al \ov R_{i,i}=\al+ E\(Z_{i}\).\label{njc1q}
   \end{equation}
Adding this up for $i=1,\ldots,n$ we get  
\begin{equation}
E\( \|Z\|_{\ell^1 } \) =\al\(\sum_{i=1}^{n} \ov R_{i,i}-n\) .\label{3.26}
  \end{equation}
 
 In the next section we give some formulas relating the $\ell^{1}$ norms of permanental random variables to the $\ell^{1}$ norms of the corresponding random variables $Z$.

\medskip  We give an alternate form of Lemma 1.3 in which the proof uses Lemma \ref{lem-1.6}.

  \begin{lemma} For any $m$ and $1\le p\le n$,  
\be E(Z_{p}^{m})=  \sum_{l=0}^{m}\sum_{(j_{0}, j_{1},\ldots, j_{l})\in J_{m}(l)}(-1)^{l+m}\al^{j_{0}}(\al+1) ^{j_{1}}\cdots (\al+l) ^{j_{l }}(R_{p,p}A_{ p,p})^{\,l}\label{1.47}
 \ee 
 where  
 \begin{equation}
J_{m}(l)=\{(j_{0}, j_{1},\ldots, j_{l})\,|\,  j_{0}+\cdots j_{l}=m; \,j_{i}\ge 1, i=0,\ldots, l-1; j_{l}\geq 0 \}.\label{1.50}
\end{equation}
 \end{lemma}

 	  \Proof  To simplify the notation we take $p=1$.  It follows from (\ref{mom.22}) that for each $m$,
 \begin{equation}
 E\(\prod_{ i=0}^{m-1}\(\al+i+Z_{1} \) \)   =(R_{1,1}A_{ 1,1})^{m}\prod_{ i=0}^{m-1}\(\al+i \).\label{10.51}
 \end{equation}
When $m=1$ this gives
 \begin{equation}
 E\(Z_{1}\)=\al  R_{1,1}A_{1,1} -\al, \label{10.51m}
 \end{equation}
 which proves (\ref{1.47}) when $m=1$. 
 
Expanding the left hand side of (\ref{10.51}) gives
\begin{equation} 
\sum_{k=0}^{m } \sum_{\stackrel{U \subseteq [0,m-1]}{|U|=m-k}}\prod_{i\in U} \(\al+i  \)E\(Z_{1}^{k}\)= (R_{1,1}A_{ 1,1})^{m} \prod_{i=0}^{m-1}\(\al+i\).\label{10.52}
\end{equation}
 We prove (\ref{1.47}) inductively. We have just seen that when $m=1$  (\ref{10.52}) holds when $E(Z_{1})$ takes the value given in (\ref{1.47}). Therefore, if we show that (\ref{10.52}) holds when $m=2$  and $E(Z_{1})$ and $E(Z_{1}^{2})$ take  the value given in (\ref{1.47}), it follows that (\ref{1.47}) gives  the correct value of $E(Z_{1}^{2})$ when $m=2$. We now assume that we have shown this up to $m-1$ and write out the left-hand side of (\ref{10.52}), replacing each $E(Z_{1}^{k})$, $k=1,\ldots,m$ by the right-hand side of (\ref{1.47}).
 Doing this we obtain terms, depending on $k$ and $U$, which, up to their sign, are of the form 
\begin{equation}
I(j_{1},\ldots, j_{m-1};l )=\prod_{ i=0}^{m-1}\(\al+i \)^{j_{i}} (R_{1,1}A_{ 1,1})  ^{\,l}\label{10.53},
\end{equation}
where $\sum_{ i=0}^{m-1}j_{i}=m$;   $0\le l\le m $;    $\,j_{i}\ge 1,  i=0,\ldots, l-1$; $   j_{l}\geq 0$ and 
$0\leq j_{i}\le 1, i=l+1,\ldots, m-1$. 

The terms  in (\ref{10.53}) may come from the term $\prod_{i\in U} \(\al+i  \)$ in (\ref{10.52}) or they may come from the expression for $E\(Z_{1}^{k}\)$ in (\ref{1.47}). Suppose that for some $ i=0,\ldots, l-1$ we have $j_{i}>1$ and  $i\in U$ and  $k=\bar k\geq l$. Consider what this term is in (\ref{10.53}). Note that we obtain the same term with a change of sign when $U$ is replaced by $ U-\{i\}$ and  $k=\bar k+1$.  The same observation holds in reverse. Furthermore, both these arguments also apply when $j_{l}> 0$.

Because of  all this cancelation, when we add up all the terms  which, up to their sign, are of the form (\ref{10.53}), and take their signs into consideration we only get non-zero contributions when
   all $\,j_{i}=1, i=0,\ldots, l-1$ and  $  j_{l}=0$.  That is, we only get non-zero contributions when 
   \begin{equation}
   \sum_{ i=0}^{l}j_{i}=l\label{10.54}
   \end{equation} 
  But recall that we have 
$\sum_{ i=0}^{m-1}j_{i}=m$ in (\ref{10.53}) so for this to hold we must have $\sum_{ i=l+1}^{m-1}j_{i}=m-l$ with $0\leq j_{i}\le 1, i=l+1,\ldots, m-1$.  This is not possible because there are only $m-l-1$ terms in this sum. Therefore we must have $l=m$ in (\ref{10.54}), which can also be written as $\sum_{ i=0}^{m-1}j_{i}=m$ because $j_{l}=j_{m}=0$.   This shows that summing all the terms on the left-hand side of (\ref{10.52}) gives the right hand side of (\ref{10.52}).   This completes the induction step and establishes (\ref{1.47}).\qed

\begin{example} {\rm It is interesting to have some examples. We have already pointed out that for any $1\le p\le n$
 \begin{equation}
E(Z _{p})=\al (A_{p,p}R_{p,p}-1). \label{1.471}
 \end{equation}
Using (\ref{1.47}) we get 
   \bea
	 \lefteqn{ E(Z_{p} ^{4})= \big[\al (\al+1) (\al+2) (\al+3)   \big](A_{p,p}R_{p,p})^{4}\label{1.474}}\\
 && \qquad-\big[\al^{2}(\al+1)  (\al+2)  +\al (\al+1) ^{2}(\al+2) \nn\\
&&  \hspace{.6 in} +\al (\al+1)  (\al+2) ^{2}+ \al (\al+1) (\al+2) (\al+3)   \big] (A_{p,p}R_{p,p})^{3}\nn\\
 &&\qquad+ \big[\al^{2}(\al+1)  (\al+2)  +\al (\al+1) ^{2}(\al+2)  +\al (\al+1)  (\al+2) ^{2} \nn\\
&& \hspace{1.1 in}+ \al^{3}(\al+1) +\al^{2}(\al+1)^{2}+  \al (\al+1) ^{3} \big] (A_{p,p}R_{p,p})^{2}\nn\\
&&\qquad-\big[ \al^{3}(\al+1) +\al^{2}(\al+1)^{2}+  \al (\al+1) ^{3}+\al^{4} \big] (A_{p,p}R_{p,p})+\al^{4} .\nn
 \eea 
Using (\ref{1.17}) we get 
\bea
	 \lefteqn{ E(Z_{p} ^{4})= \big[\al (\al+1) (\al+2) (\al+3)   \big](A_{p,p}R_{p,p})^{3}(A_{p,p}R_{p,p}-1)\label{1.474}}\\
 && \qquad-\big[\al^{2}(\al+1)  (\al+2)  +\al (\al+1) ^{2}(\al+2)\nn\\
 &&\hspace{1.3in}+\al (\al+1) (\al+2)^{2}] (A_{p,p}R_{p,p})^{2}(A_{p,p}R_{p,p}-1)\nn\\
&&\qquad+\big[ \al^{3}(\al+1) +\al^{2}(\al+1)^{2}+  \al (\al+1) ^{3} \big]   (A_{p,p}R_{p,p}) (A_{p,p}R_{p,p}-1)\nn\\
&&\qquad-   \al^{4}(A_{p,p}R_{p,p} -1) .\nn
 \eea

 }\end{example}

 	 \section{Some formulas for the $\ell^{1}$ norm of permanental random variables}

	 We can use (\ref{9.12}) and properties of independent gamma random variables to  obtain some interesting formulas for functions of permanental random variables.    
Taking advantage of the infinitely divisibility of the components of $Y$,  defined in (\ref{9.10w}), we see that
 \be
    \|Y\|_{\ell^1 }\stl \xi_{  \al n+\|Z\|_{\ell^1 },1}.\label{1.15}
  \ee  
  Note that 
  \begin{equation}
   P(\|Z\|_{\ell^1 } =j)= \sum_ {\stackrel{{\bold k}=(k_{1},\ldots, k_{n}) } {| {\bf k}|=j } }   P\(Z =(k_{1},\ldots, k_{n})\),\label{1.16}
   \end{equation}
   where $|{\bf k}|=\sum_{i=1}^{n}k_{i}$ .
  The following lemma is an immediate consequence of (\ref{1.15}):
    \begin{lemma}  \label{lem-2.1}
Let $\Phi$ be a positive real valued function. Then  
\be  
 E(\Phi(\|Y\|_{\ell^1 }))=  E(\Phi(\xi_{n\al +\|Z\|_{\ell^1 },1}))\label{3.21q} .
\ee 

\end{lemma}
 
\begin{example} {\rm 

It follows from  (\ref{3.21q}) and (\ref{1.13}) that for any $p>0$,
\begin{equation}
    E( \|Y\|^{p}_{1} )= E\(\frac{\Ga\(\|Z\|_{\ell^1 }+n\al+p\)}{\Ga\(\|Z\|_{\ell^1 } +n\al\)}\).\label{3.25qq}
  \end{equation} 

\medskip	Clearly, $  E( \|Y\|_{\ell^1 } )= \sum_{i=1}^{n} a_{i}E(X_{i})$ and $E(X_{i})=\al R_{i,i}$. Therefore, (\ref{3.26})  follows from  (\ref{3.25qq}) with $p=1$.
 
 }\end{example}

 Using these results we get a formula for the expectation of the $\ell^{2}$ norm of certain $n$-dimensional Gaussian random variables.

\begin{corollary} Let $\eta=(\eta_{1},\ldots,\eta_{n})$ be a mean zero Gaussian random variable with covariance matrix $R$. Assume that $A=R^{-1}$ exists and is an $M$-matrix. Let $\{a_{i}\}_{i=1}^{n}$ denote the diagonal elements of $A$.  Set ${\bf a^{1/2}\eta}=(a^{1/2}_{1} \,\eta_{1},\ldots,a^{1/2}_{n} \,\eta_{n})$. Then
  \begin{equation}
 \Big\|\frac{{\bf a^{1/2}\eta}}{\sqrt 2}\Big\|_{\ell^{2}}^{2}\stl   \xi_{  n/2  +  \| Z\|_{\ell^1 },1}\,.\label{11.12}
  \end{equation}
  and
   \be 
 E\( \Big\|\frac{{\bf a^{1/2}\eta }} {\sqrt 2}\Big\| _{2}\)= \label{3.21qws}     E\(\frac{\Ga\(\|Z\|_{\ell^1 }+(n+1)/2\)}{\Ga\( \|Z\|_{\ell^1 } +n/2\)}\).
\ee 
 \end{corollary}  

 \Proof  The statement in (\ref{11.12}) is simply (\ref{1.15}) with $Y=(a^{1/2}_{1} \,\eta_{1},\ldots,a^{1/2}_{n} \,\eta_{n})$ and $\al=1/2$. The statement in (\ref{3.21qws}) is simply (\ref{3.25qq}) with $p=1/2$. \qed

        \section{Symmetrizing $M$-matrices \label{sec-11}}
        
       It follows from \cite[p. 135, $G_{20}$; see also p. 150, $E_{11}$]{BP} that a symmetric  $M$-matrix is positive definite.     
Therefore when the   $M$-matrix $A$ is symmetric and $\al=1/2$,   (\ref{int.1}), with $R=A^{-1}$,  is the Laplace transform of a vector with components that are the squares of the components of a Gaussian vector.  For this reason we think of symmetric $M$-matrices as being special. Therefore, given an   $M$-matrix, we ask ourselves how does the permanental vector it defines compare with the permanental vector  defined by a symmetrized version of the   $M$-matrix. 
        
 When the $M$-matrix is symmetric and $\al=1/2$, (\ref{int.1}), with $R=M^{-1}$,  is the Laplace transform of a vector with components that are the squares of the components of a Gaussian vector.  For this reason we think of symmetric $M$-matrices as being special. Therefore, given an $M$-matrix, we ask ourselves how does the permanental vector it defines compare with the permanental vector  defined by a symmetrized version of the $M$-matrix. 
 
 For a positive $n\times n$ matrix $C$ with entries $\{c_{i,j}\}$ we define $S(C)$ to be the $n\times n$ matrix with entries $\{(c_{i,j}c_{j,i)})^{1/2} \}$.  When $A$ is an $n\times n$ non-singular  M-matrix of the form  $A=D -B$, as in (\ref{vj.73}), we set $A_{sym}=D -S(B)$.   We consider the relationship of permamental vectors determined by $A$ and $A_{sym}$, i.e. by $R=A^{-1}$ and  $R_{sym}:=A^{-1}_{sym}$ as in (\ref{int.1}).  In Remark \ref{rem-3.1} we explain how this  can be used 
 in the study of sample path properties of permanental processes.

 \begin{lemma} Let $A$ be a non-singular $M$-matrix, with diagonal elements $\{a_{i}\}$, $i=1,\ldots,n$. Then
 \begin{equation}
   |\ov A|\le 1.
   \end{equation}
 \end{lemma}
 \Proof This follows from (\ref{1.9}) since $B({\bf 0 })=1$.\qed

 The series expansion in (\ref{1.11}) gives the  following relationships between two non-singular $M$-matrices $A$ and $A'$ subject to certain regularity conditions.   
 
	\begin{lemma}\label{lem-9.2}  Let $A$ and $A'$ be $ n\times n$ non singular $M$-matrices and define $\ov A$ and ${\ov A}'$ as in (\ref{vj.73q}). Assume that ${\ov B}'\ge \ov B$.
Then
 
\begin{equation}
  {\ov A}'\le    \ov  A  \qquad\mbox{and} \qquad   ( \ov A )^{-1}\leq ( {\ov A}')^{-1}.\label{10.16q}
   \end{equation} 
     \end{lemma}

\Proof   The   first inequality  in (\ref{10.16q}) follows immediately from (\ref{vj.73q}).

To obtain the second  statement in (\ref{10.16q}) we write $A=  D_{A}(I-D_{A}^{-1}B)$, so that by \cite[Chapter 7, Theorem 5.2]{BP}  
\begin{equation}
A^{-1}D_{A}=(I-\ov B)^{-1} =\sum_{j=0}^{\ff}( \ov B)^{j}\quad\mbox{and} \quad A'^{-1}D_{A}=(I-\ov B')^{-1} =\sum_{j=0}^{\ff}( \ov B')^{j} \label{vj.76}
\end{equation}
 both converge.  Therefore $A^{-1}D_{A}\le  A'^{-1}D_{A'}$ which is the same as  the second  inequality  in (\ref{10.16q}).\qed

\begin{lemma}\label{lem-11.2} When $A$ is an $n\times n$ non-singular  M-matrix,    $A_{sym}$, and $\ov A_{sym}$ are  $n\times n$ non-singular  M-matrices and
\begin{equation}
   |A_{sym}|\ge |A|\quad\mbox{and }\quad |\ov A_{sym}|\ge |\ov A|.\label{sym}
   \end{equation}
 \end{lemma}  

\Proof  We prove this for $A$ and $A_{sym}$. Given this it is obvious that the lemma also holds for $\ov A$ and $\ov A_{sym}$.

It follows from \cite[ p. 136,  $H_{25}$]{BP} that we can find a positive diagonal matrix $E=\mbox{diag }(e_{1},\ldots, e_{n})$ such that  
\begin{equation}
EAE^{-1}+  E^{-1}A^{t}E\label{9.50p}
\end{equation}
is strictly positive definite.  We use this to  show that $A_{sym}$ is strictly positive definite.  

  We write $A=D -B$ as in (\ref{vj.73}).  For any $x=(x_{1}, \ldots,x_{n})$, by definition,
 \bea 
\sum_{i,j=1}^{n}(A_{sym})_{i,j}x_{i}x_{j}&= &\sum_{i=1}^{n}a_{i}x^{2}_{i}- \sum_{i,j=1}^{n}(b_{i,j}b_{j,i})^{1/2}x_{i}x_{j}\label{9.50}\\
&\geq & \sum_{i=1}^{n}a_{i}x^{2}_{i}- \sum_{i,j=1}^{n}(b_{i,j}b_{j,i})^{1/2}|x_{i}|\,\,|x_{j}|\nn\\
&= & \sum_{i=1}^{n}a_{i}x^{2}_{i}- \sum_{i,j=1}^{n}(e_{i}b_{i,j}e_{j}^{-1} e_{j}b_{j,i}e_{i}^{-1})^{1/2}|x_{i}|\,\,|x_{j}|,\nn
\eea
where, the first equality   uses the facts that $B\ge 0$ and has $b_{i,i}=0$, $1\le i\le n$. 

Using   the   inequality between the geometric mean and arithmetic mean of numbers we see that the last line of (\ref{9.50})
\bea
&\geq & \sum_{i=1}^{n}a_{i}x^{2}_{i}- {1 \over 2}\sum_{i,j=1}^{n}(e_{i}b_{i,j}e_{j}^{-1} +e_{j}b_{j,i}e_{i}^{-1})|x_{i}|\,\,|x_{j}| \\
&= & \sum_{i=1}^{n}a_{i}x^{2}_{i}- {1 \over 2}\sum_{i,j=1}^{n}(e_{i}b_{i,j}e^{-1}_{j} +e_{i}^{-1}b^{t}_{i,j}e_{j})|x_{i}|\,\,|x_{j}|\nn\\
& =& {1 \over 2}\sum_{i,j=1}^{n}((EAE^{-1})_{i,j}+(E^{-1}A^{t}E)_{i,j})|x_{i}|\,\,|x_{j}|> 0,\nn
\eea
by (\ref{9.50p}).  Therefore $A_{sym}$
is  strictly positive definite and by definition, $A_{sym}$   has non-positive off diagonal elements.        Since the eigenvalues of  $A_{sym}$ are real and strictly positive we see by  \cite[p. 135, $G_{20}$]{BP} that $A_{sym}$ is a non-singular  M-matrix.    
 
  To get (\ref{sym}) we note that by (\ref{1.9}) 
  \begin{equation}
| A|^{\al} \sum_{{\bold k}=(k_{1},\ldots, k_{n}) }  {| B (\bold k)|_{\al}\over  \prod_{i=1}^{n} a_{i}^{k} k_{i}! }= | A_{sym}|^{\al} \sum_{{\bold k}=(k_{1},\ldots, k_{n}) }  {| S(B) (\bold k)|_{\al}\over  \prod_{i=1}^{n} a_{i}^{k} k_{i}! }  . \label{10.13aa}
  \end{equation}
Using (\ref{9.21}) in the next lemma, we get  (\ref{sym}).\qed

\begin{lemma}\label{lem-9.3} Let $C$ be a positive $n\times n$ matrix. Then 
\be
|S(C)|_{\al}\le |C|_{\al}\quad\mbox{and}\quad|S(C)({\bold k})|_{\al}\le |C({\bold k})|_{\al}\label{9.21}.
\ee
\end{lemma}

 \Proof Consider two terms on the right-hand side of  (\ref{vj.71}) for  $|C|_{\al}$, 
\be
\al^{c(\pi)} c_{1,\pi(1)} 	c_{2,\pi(2)}\cdots  c_{n,\pi(n)}\ee
and
\be 
 \al^{c(\pi^{-1})} c_{ 1,\pi^{-1}(1) } c_{2, \pi^{-1}(2) }\cdots  c_{n, \pi^{-1}(n) }
\label{9.22u}  =  \al^{c(\pi)} c_{\pi(1),1}c_{\pi(2),2}\cdots  c_{\pi(n),n}	
\ee  The sum of these terms is
\begin{equation}
  \al^{c(\pi)}\( c_{1,\pi(1)}c_{2,\pi(2)}\cdots  c_{n,\pi(n)} 	  +  c_{\pi(1),1}c_{\pi(2),2}\cdots  c_{\pi(n),n}\).\label{9.22}
  \end{equation}
 The corresponding sum of these terms for $|S(C)|_{\al}$ is
\begin{equation}
  \al^{c(\pi)} 2 (c_{1,\pi(1)} 	c_{2,\pi(1)}\cdots  c_{n,\pi(n)}  c_{\pi(1),1}c_{\pi(2),2}\cdots  c_{\pi(n),n})^{1/2}.\label{9.23}
  \end{equation}  
 Considering the inequality between the geometric mean and arithmetic mean of numbers we see that the term in (\ref{9.22}) is greater than or equal to  the term in (\ref{9.23}). The same inequality holds for all the other terms   on the right-hand side of  (\ref{vj.71}). Therefore we have the first inequality in (\ref{9.21}).  A similar analysis gives the second inequality. \qed

 \bt\label{theo-9.3} Let $X$ and $\wt X $ be permanental vectors  determined by   $A$  and $ A_{sym}$ and $f$ be a positive function on $R^{n}$.  Then
 \begin{equation}  E(f(X))\ge  \frac{|\ov A|^{\al}}{|  \ov A_{sym}|^{\al}}   E(f(\wt X)).\label{1.28q}
  \end{equation}
 \et
 	

\Proof     Using Lemma \ref{lem-9.3} and (\ref{9.12}) we  have
\bea
\lefteqn{ E(f(  X))\label{9.12a}}\\
  &&=   |\ov A | \sum_{{\bold k}=(k_{1},\ldots, k_{n}) }  {|\ov B (\bold k)|_{\al}\over  \prod_{i=1}^{n}  k_{i}! }E\( f\(\xi_{\al+k_{1},a_{1}} ,\ldots,\xi_{\al+k_{n},a_{n}} \)\)\nn\\
  &&\geq  |\ov A |  \sum_{{\bold k}=(k_{1},\ldots, k_{n}) }  {|\ov B_{sym}(\bold k)|_{\al}\over  \prod_{i=1}^{n}  k_{i}! }E\( f\(\xi_{\al+k_{1},a_{1}} ,\ldots,\xi_{\al+k_{n},a_{n}} \)\)\nn\\
  &&=\frac{|\ov A|^{\al}}{|\ov A_{sym}|^{\al}}E(f(\wt X )).\nn 
  \eea 
\qed

  This leads to an interesting two sided inequality.

\begin{corollary}\label{cor-7.1} Let $X$ and $\wt X $ be permanental vectors  determined by   $A$  and $ A_{sym}$. Then for all functions $g$ of $X$ and $\wt X$ and sets $\BB$ in the range of $g$
\bea
 \frac{|A|^{\al}}{|  A_{sym}|^{\al}}  P\(g( \wt X)\in \BB\) &\le &P\( g(   X)\in\BB\)\label{3.16w}\\
 & \le & \(1- \frac{|A|^{\al}}{|  A_{sym}|^{\al}} \) +\frac{|A|^{\al}}{|  A_{sym}|^{\al}}P\(g( \wt X)\in\BB\).\nn
 \eea
 \end{corollary}

  \Proof   The first  inequality follows by taking  $f(X)=I_{g( X)\in \BB}(\cd)$ in (\ref{1.28q}) and, similarly, the second inequality follows by taking $f(X)= 	I_{g( X)\in \BB^{c}}(\cd)$ in (\ref{1.28q}). \qed
  
   \begin{corollary}\label{cor-3.2} Under the hypotheses of Corollary \ref{cor-7.1}\begin{equation}
P\(g( X)\in \BB\) =1\implies  P\(g(\wt  X)\in \BB\)=1 .   \label{3.16}
   \end{equation}\end{corollary}
   
   \Proof  It follows from the first inequality in (\ref{3.16w}) that  
   \begin{equation}
P\(g( X)\in \BB^{c}\) =0\implies  P\(g(\wt  X)\in \BB^{c}\)=0.   \label{3.16s}
   \end{equation}
  We get  (\ref{3.16}) by taking complements. \qed
  
   \begin{remark} {\rm A useful application of Corollaries \ref{cor-7.1} and \ref{cor-3.2} is to take $g(X)=\| X\|_{\ff}$.

 }\end{remark}

\begin{remark}\label{rem-3.1} {\rm  When $\{R(s,t),s,t\in S\}$ is the potential density of a transient Markov process with state space $S$, for all $(s_{1},\ldots,s_{n})$ in $S$, the matrix   $\{R(s_{i},s_{j})\}_{i,j=1}^{n}$ is invertible and  its inverse  $A(s_{1},\ldots,s_{n})$  is a non-singular $M$-matrix. For all $(s_{1},\ldots,s_{n})$ in $S$ consider $A_{sym}(s_{1},\ldots,s_{n})$.  
 If  
\begin{equation}
 \inf_{\forall t_{1},\ldots,t_{n}, \forall n}   \frac{|A(t_{1},\ldots,t_{n} )| }{|  A_{sym}(t_{1},\ldots,t_{n} )| }>0
   \end{equation}
it follows from Corollary \ref{cor-7.1} that $\sup_{t\in T}X_{t}<\ff$ almost surely if and only if $\sup_{t\in T}\wt X_{t}\newline <\ff$ almost surely. Here we also use the fact that these are tail events.
 }\end{remark}

\bibliographystyle{amsplain}

      \def\wh{\widehat}
\def\ol{\overline}

\end{document}